\documentclass{article}
\usepackage[utf8]{inputenc}
\usepackage{amsmath}
\usepackage{amssymb}
\title{On the Frame-Stewart Conjecture}
\author{Youngjin Bae }
\newtheorem{defn}{Definition}[section]
\newtheorem{ex}{Example}[section]

\newtheorem{thm}{Theorem}[section]
\newtheorem{con}{Conjecture}[section]
\newtheorem{lemma}{Lemma}[section]
\date{}

\begin{document}

\maketitle

\begin{abstract}

The Frame-Stewart conjecture states the least number of moves to solve a generalized Tower of Hanoi problem, of n disks and p pegs. 
In this paper, we prove several weaker versions of the Frame-Stewart conjecture.

\end{abstract}
\vspace{7mm}
\tableofcontents
\newpage
\section{Introduction}

\vspace{5mm} 

The generalized Tower of Hanoi problem can be formally stated as following.

\begin{defn}
Let $n$ and $p\geq3$ be natural numbers. Than a generalized Tower of Hanoi problem is a problem of moving n ordered disks (we may number those disks from 1 to n) from an initial peg to another one, satisfying following conditions:

1. No larger disk can be on top of a smaller one

2. A disk can be moved from one peg to another peg only when no other 

disks are on top of it. 

\end{defn}

We simply call a generalized Tower of Hanoi problem with $n$ disks and $p$ pegs as $(n,p)-problem$. 

\begin{defn}

For $n$ and $p\geq3$, $M(n,p)$ is the least number of moves needed to solve $(n,p)-problem$.
\end{defn}

\begin{thm}[The original Tower of Hanoi problem]
For natural number $n$, $$M(n,3)=2^n-1$$
\end{thm}

\begin{thm}[A.A.K.Majumdar]
For $n,p$, there exist an unique natural number $r$ satisfying 

$$ \binom{p+r-3}{p-2} \leq n <  \binom{p+r-2}{p-2} $$

and 

$$M(n,p) \leq \sum_{t=0}^{r-1}{2^{t} \binom{p+t-3}{p-3}} +2^{r}(n-\binom{p+r-3}{p-2})$$
\end{thm}

\textbf{Proof.} (See [2], for example) The existence of $r$ is clear. Define $ K(n,p)=\sum_{t=0}^{r-1}{2^{t} \binom{p+t-3}{p-3}} +2^{r}(n-\binom{p+r-3}{p-2}) $. We will prove the inequality by explicitly showing that it is possible to solve the $(n,p)-problem$ with exactly $K(n,p)$ times of move. (1)

We use induction on $p$ and then on $n$. 
First, for $p=3$, we have $M(n,3)=2^n-1=K(n,3)$ for all $n$. 
Assume that it is possible to solve the $(n,p)-problem$ with $K(n,p)$ times of move for $3 \leq p \leq q-1$. For $p=q$, we use induction on $n$. For $n=1$, $K(1,q)=1$ and it is indeed possible to move a single disk with 1 move. 
Assume that (1) holds for $n \leq m-1$. For $n=m$, let $m=\binom{q+r-3}{q-2}+\alpha$ where $0 \leq \alpha < \binom{q+r-3}{q-3} $. Since $\binom{q+r-3}{q-3}=\binom{q+r-4}{q-4}+\binom{q+r-4}{q-3}$, there are natural numbers $\beta, \gamma$ such that $\alpha=\beta+\gamma$ and $\beta < \binom{q+r-4}{q-3}$ and $\gamma < \binom{q+r-4}{q-4}$. 
We call the peg on where every disks are at the beginning as initial peg ($I$ for short) and the peg on where every disks are at the end as final peg ($F$ for short). Also, since $p \geq 3$, we can pick a peg different from $I$ and $F$ and call it middle peg ($M$ for short). Note that $m=\binom{q+r-3}{q-2}+\alpha = (\binom{q+r-4}{q-2}+\beta) +(\binom{q+r-4}{q-3}+\gamma)$. Define $k:=\binom{q+r-4}{q-2}+\beta$ and we have $m-k=\binom{q+r-4}{q-3}+\gamma$. Then, we move $m$ pegs through the following process:

\vspace{5mm}
1. Move disks 1 to k from $I$ to $M$ with $K(k,q)$ moves. 

2. Move disks k+1 to m from $I$ to $F$ with $K(m-k,q-1)$ moves. (Note that we do not use the peg $M$ here.)

3. Move disks 1 to k from $M$ to $F$ with $K(k,q)$ moves.

\vspace{5mm}

So far, we have moved the $m$ disks with $2K(k,q)+K(m-k,q-1)$ moves. Now it is enough to check that $$K(m,q)=2K(k,q)+K(m-k,q-1)$$. 
\linebreak
This can be shown by calculation:

\vspace{5mm}
We have $$K(k,q)=\sum_{t=0}^{r-2} 2^{t}\binom{q+t-3}{q-3}+2^{r-1}\beta$$ and $$K(m-k,q-1)=\sum_{t=0}^{r-1} 2^{t}\binom{q+t-4}{q-4}+2^{r}\gamma$$.
\linebreak
Thus, $$2K(k,q)+K(m-k,q-1)=\sum_{t=0}^{r-2} 2^{t+1}\binom{q+t-3}{q-3}+2^{r-1}+\sum_{t=0}^{r-1} 2^{t}\binom{q+t-4}{q-4}+2^{r}(\beta+\gamma)$$ 
$$=\sum_{t=0}^{r-1} 2^{t}(\binom{q+t-4}{q-3}+\binom{q+t-4}{q-4})+2^{r}\alpha=K(m,q)$$

\vspace{5mm}
Which finishes the proof. Note that the proof works for every possible $\beta$ and $\gamma$ satisfying the conditions, which implies that the minimal solution might not be unique.

\section{The Frame-Stewart Conjecture}

The Frame-Stewart Conjecture states that the DP-algorithm in the proof of previous theorem is actually optimal and thus $M(n,p)=K(n,p)$.
\begin{con}[Frame-Stewart Conjecture]
For $n,p$, $M(n,p)=K(n,p)$ holds. 
\end{con}

The conjecture indeed holds for $p=3$.

\section{Preliminary Facts}

For natural number $x$, we define $\bar{x} :=\{1,2,..,x\}$.

\begin{defn}
Given $n,p$, a \textbf{state} of the $(n,p)-problem$ ($(n,p)-state$ in short) is $n$ disks being allocated on $p$ pegs. Formally, a state is equivalent to a function $f:\bar{n} \to \bar{p}$ We define the set of all states of the $(n,p)-problem$ as 

$$X(n,p):=\{f:\bar{n} \to \bar{p}\}$$
\end{defn}

\begin{defn}
Given $n,p$ and two states $f,g$ of the $(n,p)-problem$, a \textbf{path} connecting $f$ and $g$ is a finite sequence of $(n,p)-states$ such that the initial term of the sequence is $f$ and the final term is $g$. 
If $P=\{P_{0}=f, P_{1}, ..., P_{k}=g\}$ is a path connecting(between) $f$ and $g$, define \textbf{length} of the path as $|P|:=k$. 
\end{defn}

\begin{defn}
Let $f,g$ be $(n,p)-states$. Define $P(f,g)$ as the set of all paths connecting $f$ and $g$. A path between $f$ and $g$ is a \textbf{shortest path} if its length is minimal among $P(f,g)$. A length 1 path is called \textbf{move} from $f$ to $g$. We formally write a shortest path between $f$ and $g$ as $f-g$. It is obvious that $|f-g|=|g-f|$ for any given $f$ and $g$. Note that shortest path between $f$ and $g$ may not unique and $f-g$ is not well-defined. Still, $|f-g|$ is well-defined. 
\end{defn}

\begin{defn}
Let $f,g$ be $(n,p)-states$ and $\psi$ be a path between $f$ and $g$. If $X$ is a subset of $\bar{n}$, we define $|\psi|_{X}$ be the number of moves of disks in $X$ while $\psi$. 
\end{defn}

\begin{ex}
$(n,p)-problem$ can be demonstrated as finding shortest path between two distinct constant states (i.e. constant function) $f$ and $g$. 
\end{ex}

We introduce a notation by Roberto Demontis and a notion of demolishing sequence. The triple $(j,i,t)$ with $1 \leq j < i \leq \infty$ and $j<t \leq \infty$, denotes that the disk $j$ moves from being on the disk $i$ to be placed on the disk $t$. We write $i=\infty$ when there was no disk under $j$ before it moves onto $t$. Similarly, we write $t=\infty$ when disk $j$ moves to an empty peg. 

\begin{defn}
A path $P$ between $f$ and $g$ is said to be \textbf{demolishing sequence} if 

1. $f$ is a constant state

2. The last move of $P$ is $(n,\infty,\infty)$

3. The move $(n,\infty,\infty)$ appears exactly once in $P$.

\end{defn}

We call the final state of a minimal demolishing sequence as \textbf{middle state}. 

\begin{defn}
Let $P$ and $Q$ be sequences satisfying $P_{|P|}=Q_{0}$. Define $P+Q$ be a sequence concatenate $P$ and $Q$. $|P+Q|=|P|+|Q|$ holds. 
\end{defn}

\begin{thm}[Roberto Demontis]
Given $f \equiv I_{f},g \equiv I_{g}$ be two distinct constant states and $S:=f-g$. Assume $f-h$ be a minimal demolishing sequence of moves. Then, $|S|=2|f-h|+1$ holds. 
\end{thm}

\textbf{Proof} Since $f$ and $g$ are two distinct constant states, there must be at least one $(n,\infty,\infty)$ move in $S$, which we will call $\psi$. Let $P$ be a subsequence of $S$ from the beginning to the last move before $\psi$. Than, $P$ is a demolishing sequence. Let $S=P+\psi+Q$. Define $P^{r}$ be a sequence which is reverse of $P$ but $I_{f}$ and $I_{g}$ are switched. If $|P|>|Q|$, we have $|S|=|P+\psi+Q|<|P+\psi+P^ {r}|$, which contradicts to the minimality of $S$. Similarly, if $|P|<|Q|$, we have $|S|=|P+\psi+Q|<|Q^{r}+\psi+Q|$, also contradiction. Thus, we have $|P|=|Q|$ and $|S|=2|Q|+1=2|f-h|+1$ since both $Q$ and $f-h$ are minimal demolishing sequences.

By the theorem above, it is enough to find minimal demolishing sequence instead of the whole $(n,p)-problem$.

\begin{thm}[Roberto Demontis]
Let $S$ be a minimal demolishing sequence of $(n,p)-problem$. Suppose that the disks have been arranged on $r \leq p-1$ stacks at the end of $S$. Let $n,n-1$ and $j_{1} < ... < j_{r-2}$ be the disks at the bottom of the $r$ stacks at the end of $S$. Then during the demolishing phase, no disk $y>j_{1}$ has arranged on the peg on which the disk $j_{1}$ will be stacked at the end of $S$.
\end{thm}

\textbf{Proof} See [1].

\section{Main Results}

\begin{defn} 
Let $\mu$ be a middle state of a solution $S$ of $(n,p)-problem$. Assume that $k<n$ be the largest disk which is not stacked on $n-1$ at $\mu$. We define $B(S)=n-k-1$ the \textbf{base} of $S$. 
\end{defn}

The above definition implies that every disks $j$ of $k+1 \leq j < k+1$ are stacked on $n-1$ at $\mu$.

\begin{defn}
For a sequence $P$ and a state $\mu \in P$, let $\mu^{+}$ be the next state of $\mu$ and $\mu^{-}$ be the state before $\mu$. 
\end{defn}

\begin{defn}
$\chi$ be a sequence of $(n,p)-problem$ and $A \subset \bar{n}$. Let $\chi|_{A}$ be a sequence of $(A,p)-problem$ such that disks of $\bar{n} \backslash A$ are removed from $\chi$. $\chi|_{A}$ is called restriction of $\chi$ on $A$.
\end{defn}

\begin{lemma}
Let $n,p$ are natural numbers satisfying $\binom{r+p-3}{p-2} \leq n < \binom{r+p-2}{p-2}$. Then, $K(n,p)-K(n-1,p)$ is either $2^r$ or $2^{r-1}$. 
\end{lemma}

\textbf{Proof} If $n=\binom{r+p-3}{p-2}$, then $n-1=\binom{r+p-3}{p-2}-1$. 

We have $$K(n,p)=\sum_{t=0}^{r-1} 2^{t}\binom{p+t-3}{p-3}$$

and

$$K(n-1,p)=\sum_{t=0}^{r-2}\binom{p+t-3}{p-3}+2^{r-1}(\binom{p+r-3}{p-2}-1-\binom{p+r-4}{p-2})$$

Thus 

$$K(n,p)-K(n-1,p)=2^{r-1}$$. 

Otherwise, $\binom{p+r-3}{p-2}<n<\binom{p+r-2}{p-2}$ holds. It is obvious that $K(n,p)-K(n-1,p)=2^{r}$.

\begin{lemma}
For a sequence $\chi$, $|\chi| \geq |\chi_{0}-\chi_{|\chi|}|$ holds. Equality holds when $\chi$ is minimal. 
\end{lemma}

\begin{thm}
For a solution $S$ of $(n,p)-problem$, if $B(S) \geq r$ where $r$ is the unique natural number satisfying $\binom{r+p-3}{p-2} \leq n < \binom{r+p-2}{p-2}$, $|S| \geq K(n,p)$. In other words, if there is a shorter solution $S$ which contradicts to the Frame-Stewart conjecture, than it must satisfy $B(S) < r$. 
\end{thm}

\textbf{Proof} Let $j$ be the initial state and $\mu$ be the middle state of $S$. Define $\nu$ be the state when the $\bar{n-1} \backslash \bar{k}$ tower completes. i.e. the state right after the last $(k+1,*,k+2)$ move between $j$ and $\mu$. 

First, in case of $|j-\nu|_{\bar{k}} \geq \frac{K(n,p)}{2}-K(B(S),3)$, we have $$M(n,p)=|S|=2|j-\mu|+1=2(|j-\mu|_{\bar{k}}+|j-\mu|_{\bar{n-1} \backslash \bar{k}})+1 \geq 2|j-\nu|_{\bar{k}}+2|j-\mu|_{\bar{n-1} \backslash \bar{k}}+1$$

By \textbf{Theorem 3.2}, through the sequence $j-\mu$, any disks in $\bar{n-1} \backslash \bar{k}$ have not placed on any other pegs than the initial peg, $\mu(n-1)$ and $\mu^{+}(n)$ where $\mu^{+}$ is the state right after $\mu$. Therefore, we have $$|j-\mu|_{\bar{n-1} \backslash \bar{k}} \geq K(B(S),3)$$

Thus, in this case, $2|j-\nu|_{\bar{k}}+2|j-\mu|_{\bar{n-1} \backslash \bar{k}}+1 \geq 2(\frac{K(n,p)}{2}-K(B(S),3))+2K(B(S),3)+1 \geq K(n,p)$ and $|S| \geq K(n,p)$. 

Otherwise, $$|j-\nu|_{\bar{k}} < \frac{K(n,p)}{2}-K(B(S),3)$$ holds. Define $T:=j-\mu$ and we have $|S|=2|T|+1$. Since $\nu$ is the state right after the $\bar{n-1} \backslash \bar{k}$ tower has completed, no disks of $\bar{n-1} \backslash \bar{k}$ has moved after $\nu$.
Let $\chi$ be a sequence such that $\chi_{0}=\nu$,  $\chi_{|\chi|}=\mu(n-1)$ and $\chi|_{\bar{k}}=\nu|_{\bar{k}}-\nu(n-1)$ (i.e. $\chi$ is a sequence beginning from $\nu$ and moving disks of $\bar{k}$ onto $n-1$ minimally, instead of end up with $\mu$.)
Note that $|\chi| \leq |j-\nu|_{\bar{k}}$ holds. This is because $\chi_{0}|_{\bar{k}}=\nu|_{\bar{k}}$, $\chi_{|\chi|}|_{\bar{k}}=j|_{\bar{k}}$ and \textbf{Lemma 4.2}. 
However, we also have the sequence $(j-\nu)+\chi$, which begins with $j$ and end up with complete $\bar{n-1}$ tower. This gives the following:
$$|(j-\nu)+\chi| \geq K(n-1,p)$$
Thus we get

\vspace{5mm}

$$|T|+|j-\nu|_{\bar{k}} \geq |j-\nu|+|\chi| \geq K(n-1,p)$$

\vspace{5mm}
\newpage
and 

$$|S|=2|T|+1 \geq 2(K(n-1,p)-|j-\nu|_{\bar{k}})+1 \geq 2K(n-1,p)-2(\frac{K(n,p)}{2}-K(B(S),3)-1)+1$$
$$=2K(n-1,p)-K(n,p)+2K(B(S),3)+3 $$
\newline
By \textbf{Lemma 4.1}, $2K(n-1,p)-K(n,p) \geq K(n,p) - 2^{r+1}$ and thus 

$$2K(n-1,p)-K(n,p)+2K(B(S),3)+3 \geq K(n,p)+2K(B(S),3)+3-2^{r+1}$$
$$=K(n,p)+2(2^B(S)-2^r)+1 \geq K(n,p)$$, 

which finishes the proof.

\vspace{5mm}
The only case left is $B(S) \leq t$ for proving the Frame-Stewart Conjecture. 

\begin{con}
For $n,p$ such that $\binom{r+p-3}{p-2} \leq n < \binom{r+p-2}{p-2}$ and $S$ be a solution of $(n,p)-problem$. $|S| \geq K(n,p)$ holds for $B(S) <r$
\end{con}

\newpage


\end{document}